\documentclass[12pt,reqno]{amsart}

\usepackage{amssymb,latexsym,bbm,graphicx,epsfig,epic,eepic,oldgerm,psfrag,graphics}
\usepackage{amsthm,amsfonts,amsmath,amscd}

\usepackage[all]{xy}  

\theoremstyle{plain}

\newtheorem*{remark*}{Remark}
\newtheorem*{remarks*}{Remarks}

\newtheorem*{example*}{Example}
\newtheorem*{examples*}{Examples}

\def\1{\:\!}
\def\2{\;\!}
\def\s{\smallskip}
\def\m{\medskip}

\def\im{\operatorname {im}}

\def\Fred{\operatorname{Fred}}

\def\trace{\operatorname{trace}}

\def\ind{\operatorname{ind}}

\def\span{\operatorname{span}}

\def\FH{\operatorname{FH}}
\def\symp{\operatorname{symp}}
\def\vol{\operatorname{vol}}

\def\SL2{\operatorname{SL_2}}

\def\ga{\alpha}
\def\gb{\beta}
\def\gg{\gamma}

\def\gf{\varphi}

\def\cm{{\mathcal M}}

\def\cp{{\mathcal P}}

\def\NN{\mathbbm{N}}

\def\QQ{\mathbbm{Q}}
\def\RR{\mathbbm{R}}

\def\ZZ{\mathbbm{Z}}

\def\pp{\partial}

\def\ni{\noindent}
\def\b{\bigskip}
\def\m{\medskip}

\def\id{\mbox{id}}

\begin{document}

\title[]{The beginnings of symplectic topology in Bochum in the early eighties}

\author{Eduard Zehnder}  
\address{
ETH Z\"urich,
Switzerland}
\email{zehnder@math.ethz.ch}

\date{\today}

\maketitle

\begin{quote}
{\footnotesize
I outline the history and the original proof of the Arnold conjecture on fixed
points of Hamiltonian maps for the special case of the torus, leading to a sketch of
the proof for general symplectic manifolds and to Floer homology.
This is the written version of my talk at the Geometric Dynamics Days~2017 
(February 3--4) at the RUB in Bochum. 
I would like to thank Felix Schlenk for improvements and for his enormous help 
in typing a barely readable manuscript.
}
\end{quote}

\smallskip
\bigskip
Peter Albers has asked me to recall the beginnings of symplectic topology
here in Bochum during the early eighties. It is history, a story dating back more than
thirty years.

This reminds me of a talk I gave at the Moscow Mathematical Seminar.
The audience had decided that my talk should be translated into Russian.
So I asked Vladimir Arnold, who had invited me:
``How do I proceed? Shall I say several sentences and then wait for the translator?''
Arnold answered immediately:
``Don't worry, the translator will always be ahead of you.''

I am afraid that for the next hour all of you will be in the same situation as the translator, 
namely always ahead of me.

\b \ni
{\bf 1980 -- Forced oscillations on $\RR^{2n}$}

\m \ni
In the late seventies and early eighties of the previous century, Herbert Amann
and I were looking for forced oscillations of time-dependent Hamiltonian equations
on the standard symplectic space~$\RR^{2n}$.
Forced oscillations are 1-periodic solutions~$x(t)$ of a Hamiltonian system
$$
\frac{d}{dt} x(t) \,=\, J \1 \nabla H(t,x(t)) \quad \mbox{satisfying} \quad x(t+1)=x(t) .
$$ 
Here $H \colon \RR \times \RR^{2n} \to \RR$ is a smooth Hamiltonian function
which is also periodic in $t$ of period~one, $H(t+1,x) = H(t,x)$.
The matrix~$J$ is the standard symplectic matrix 
$$
J \,=\ 
\left( {\begin{array}{cc}
   0  & 1 \\
   -1 & 0 \\
        \end{array} } \right)
\quad \mbox{satisfying}  \quad J^2  = - 1.
$$
The flow $\{ \varphi^t \}$ consists of symplectic diffeomorphisms of~$\RR^{2n}$ and is 
defined by the solutions of the Cauchy initial value problems
$$
\frac{d}{dt} \varphi^t(x) \,=\, J \1 \nabla H(t,\varphi^t(x)) \quad \mbox{and} \quad \varphi^0 (x) = x.
$$
If $t=1$, the diffeomorphism $\varphi^1 = \varphi_H^1$ is called a Hamiltonian map, 
and the forced oscillations correspond to the fixed points of~$\varphi^1$, namely 
$\varphi^1(x(0)) = x(1) = x(0)$.

From physics we know that the forced oscillations are characterized as the critical points of the so-called 
action functional $f \colon C^\infty (S^1, \RR^{2n}) \to \RR$,
defined on the space of smooth parametrized loops on~$\RR^{2n}$ by
$$
f(u) \,=\, \int_0^1 \Bigl\{ \tfrac 12 \langle -J \dot u (t), u(t) \rangle - H (t,u(t)) \Bigr\}\, dt 
$$
Indeed, taking the derivative in the direction of the smooth loop~$v$ we find 
$$
df(u) \cdot v \,=\, \frac{d}{d \varepsilon} f(u + \varepsilon v) \big|_{\varepsilon =0} 
\,=\, \int_0^1 \bigl\langle -J \dot u(t) - \nabla H(t,u(t)) , v(t) \bigr\rangle \, dt .
$$
We see that the $L^2$-gradient of $f$ is the loop 
$$
\nabla f(u) (t) \,=\, - J \dot u (t) - \nabla H (t,u(t)) ,
$$
and we conclude that the variational principle singles out the distinguished loops in the loop space
which are the forced oscillations.

In sharp contrast to the variational principle for closed geodesics in Riemannian geometry, 
the action functional -- due to the symplectic structure -- is in general neither bounded from above nor from below.
Indeed, take the special loops
$$
u_k(t) \,=\, e^{2\pi k J t} \mathbf{e} \,=\, \cos (2\pi kt) \2 \mathbf{e} + \sin (2\pi kt) J \mathbf{e} 
$$
for a unit vector $\mathbf{e} \in \RR^{2n}$.
Then $\|u_k\|_{L^2} =1$ and 
$$
\int_0^1 \tfrac 12 \bigl\langle -J \dot u_k(t), u_k(t) \bigr\rangle \, dt = \pi k
$$
which diverges to $\pm \infty$ as $k \to \pm \infty$.
Therefore, the standard approaches to find critical points, 
like the direct method of the calculus of variations and
Morse theory in Banach spaces, do not apply to the action functional, 
which does not seem suitable for existence proofs.

The variational principle $df (u) = 0$ 
was used in physics lectures 
as the principle of stationary action in the Hamiltonian formalism
more for philosophical reasons and not for existence proofs. 
Only towards the end of the seventies Paul Rabinowitz constructed
sophisticated critical point techniques of mountain pass type in Banach spaces
and demonstrated -- very much against the advice of his advisors at the Courant Institute --
that the principle can be used very efficiently for existence proofs, see~\cite{Rab78}.
Afterwards, many people used these techniques to prove the existence of forced oscillations
in~$\RR^{2n}$,
but under quite artificial assumptions on the non-linearity at infinity.

With H.\ Amann, we therefore looked for forced oscillations for the also artificial, but this time
asymptotically linear Hamiltonian systems.
We applied H.\ Amann's so-called saddle point reduction from~\cite{Amann79, AmannZehnder80}
to the action functional to obtain a variational functional on a finite dimensional subspace, 
to which we then applied Charles Conley's index theory.
To learn more about his theory, I invited my friend Charley to Bochum to lecture about dynamical systems.

\begin{figure}[ht]
 \begin{center}
  \includegraphics[width=4.5cm]{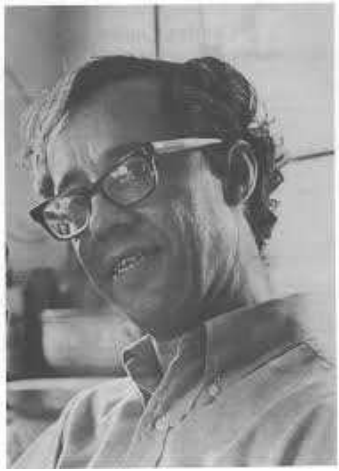}
 \end{center}
 \caption{Charles Conley}
\end{figure}

\m
Conley was an extremely unconventional, sometimes eccentric man, who created his own mathematics.
You never quite knew whether a statement was a joke or meant seriously. 
He could quote Mark Twain in every situation. He had learned Twain's work by heart while serving four and a half years
in the air force, guarding a junk yard in England. 

Charles was never very self-confident. Studying mathematics at MIT, he got the impression that all the other students
were much more clever. 
Therefore, he did not show up at the lectures any more. 
J\"urgen Moser, at that time professor at MIT, called him at home and said just one sentence:
``Mr.\ Conley, your seminar talk is scheduled for Wednesday, 5~p.m.''
That such a great man actually called him personally and even thought that he had something interesting to say,
was such a boost for Charley that he continued his studies.

As soon as Charley arrived in Bochum, we started working on a Morse theory for forced
oscillations. We first replaced the Morse index, which does not exist, 
by an index associated intrinsically to a non-degenerate forced oscillation and denoted by
the Greek letter 
$$
\mu \colon \cp_H = \{ \mbox{ forced oscillations }\} \to \ZZ ,
$$
see \cite{CoZe84}.
This integer, that is defined by the linearized Hamiltonian flow along the forced oscillation, 
describes the mean winding of nearby solutions around the forced oscillation during one period.

One day, John Mather visited us, coming from Paris where he had talked to Michel Herman.

\begin{figure}[ht]
 \begin{center}
  \includegraphics[width=6.6cm]{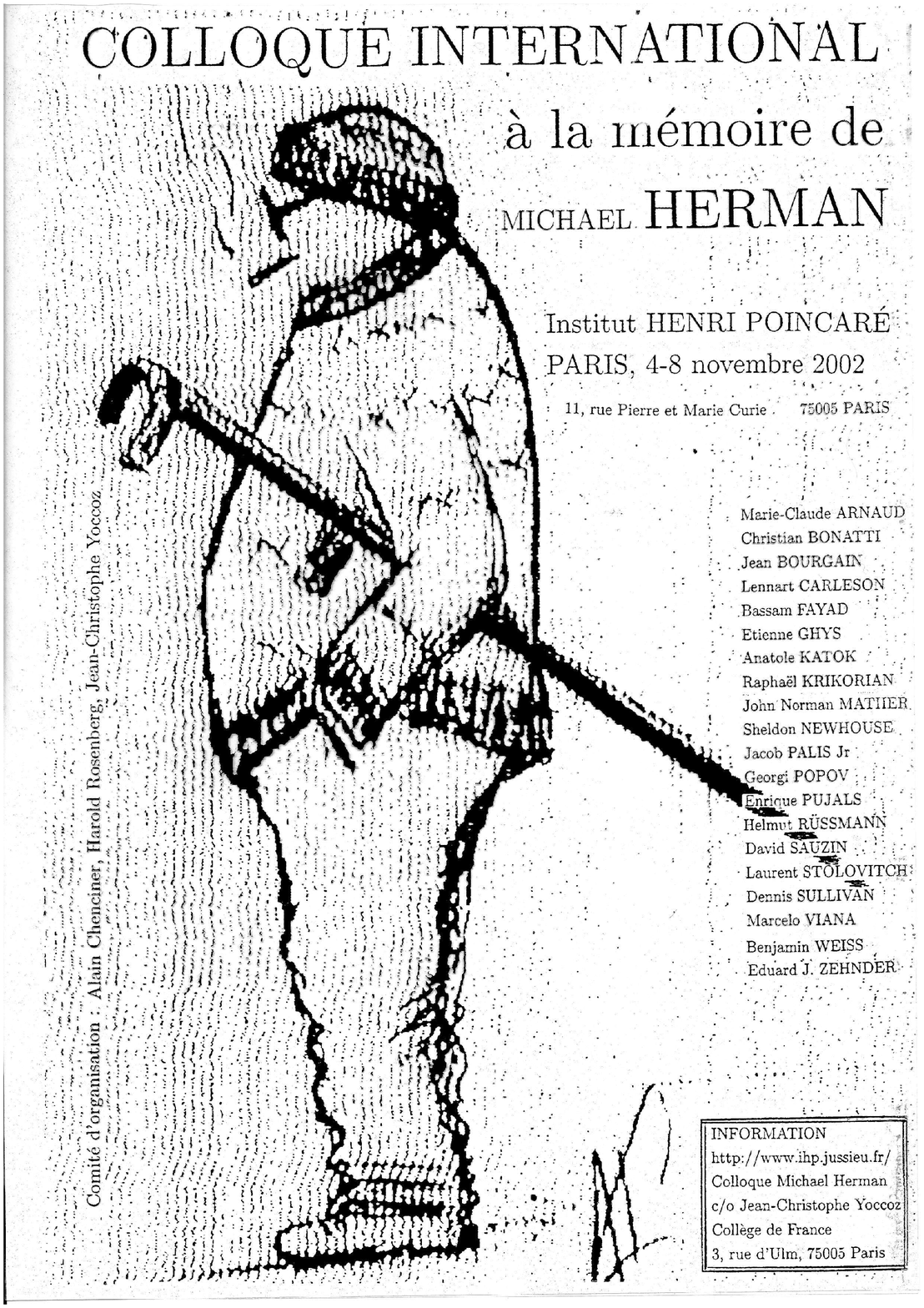}
 \end{center}
\end{figure}

\ni
Here is a picture of M.\ Herman drawn by Marie-Jo L\'ecuyer, his secretary, 
like a figure out of a Dickens novel.
Obviously, he appreciated good food, the crutch always under his arm,
and a smoking Gauloise between his lips, a sign that he was thinking.

Michel Herman was the sharpest mind in classical dynamical systems I have ever met, so far.
Of course, we inquired immediately what Michel was working on. 
John Mather told us that Michel together with his group of brilliant students tried very hard,
and in the end not successfully, to decipher an unreadable and un-understandable manuscript
of an unknown Russian mathematician named Eliashberg, 
about the Arnold conjecture for surfaces, in particular the 2-torus.
So, to be polite Charley asked him what this conjecture was.

\b \ni
{\bf V.\ Arnold's conjecture for $T^2$}

\s \ni
{\it A Hamiltonian diffeomorphism $\varphi$ of the standard torus~$T^2$
possesses at least as many fixed points as a smooth function on~$T^2$ has 
critical points, namely

$\geq 3 = \mbox{cup length}\2 (T^2)$ (by Ljusternik--Schnirelman), and

$\geq 4 =$ sum of the Betti-numbers (by Morse-theory)
if all the fixed 

points are non-degenerate.
}

\begin{figure}[ht]
 \begin{center}
  \psfrag{c}{$T^2 = \RR^2/\ZZ^2$}
 \leavevmode\epsfbox{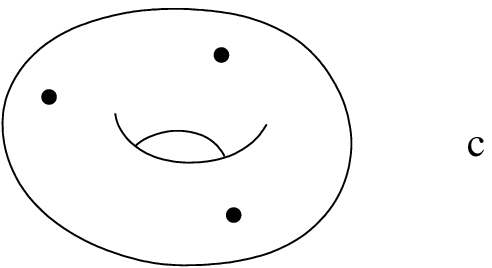}
 \end{center}
\end{figure}
%
%

\m \ni
Arnold made this conjecture first in 1966, see \cite{Arn66}
and also~\cite[p.\ 284]{Arn04}.

\b
Charley and I looked at each other and told John that we can give him
a proof for the torus in all dimensions~$2n$ in a few days.
When we showed him the proof he said:
``You are very fast in proving old conjectures, how about proving another old conjecture, 
the Riemann Hypothesis?''
Again, we asked what that conjecture was,
but then, obviously, did not touch it.

\m
Before I sketch our proof in which we used just what we had at our fingertips at
the time, I would like to quickly recall what was known to us about fixed points of Hamiltonian or symplectic maps and where the Arnold conjecture comes from.

Known to us at the time was the Lefschetz topological fixed point theory.
A continuous map $f \colon M \to M$ of a compact manifold has a fixed point
if its Lefschetz number
$$
L(f) \,=\, \sum_{k \geq 0} (-1)^k \trace \2 (f_* | H_k(M;\QQ)) \neq 0 .
$$
If $f$ is homotopic to the identity, $f \sim \id$, then 
$L(f) = L(\id) = \chi (M)$, the Euler characteristic of~$M$.
Hence, a continuous map $f \colon S^2 \to S^2$ homotopic to the identity
possesses always a fixed point.
It might have only one, degenerate fixed point:

\begin{figure}[ht]
 \begin{center}
  \psfrag{c}{$\chi(S^2)=2$}
 \leavevmode\epsfbox{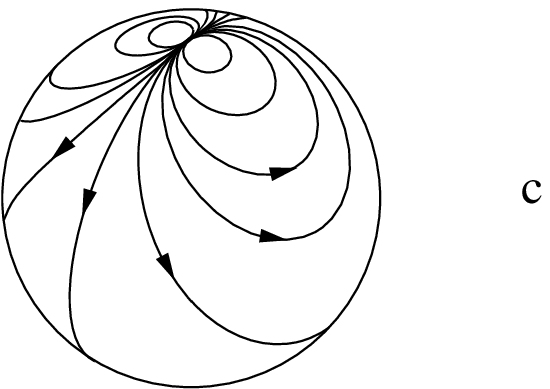}
 \end{center}
\end{figure}

\ni
However, if we add more structure and require, in addition to $f \sim \id$,
that $f$ is area preserving, then it possesses at least two fixed points
by Brouwer's translation theorem in~$\RR^2$.
If $f \colon S^2 \to S^2$ is a symplectic diffeomorphism, $f^*\omega = \omega$,
then $f \sim \id$ by a theorem of Heinz Hopf, and hence $f$ possesses at least two fixed points.
The exceptional surface is the torus $T^2 = \RR^2 / \ZZ^2$,
whose Euler characteristic vanishes, $\chi (T^2) =0$.

\begin{figure}[ht]
 \begin{center}
  \psfrag{0}{$0$}
  \psfrag{1}{$1$}
  \psfrag{11}{$(1,1)$}
 \leavevmode\epsfbox{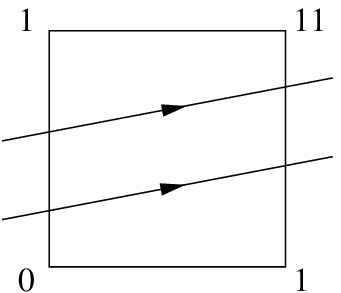}
 \end{center}
\end{figure}

\ni
Non-trivial translations on $T^2$, for example, are area preserving
and have no fixed points,
and irrational translations not even have periodic points.

\m
The most celebrated two-dimensional fixed point theorem is the

\m \ni
{\bf Poincar\'e--Birkhoff fixed point theorem (1912/13)}

\s \ni
{\it 
An area and orientation preserving homeomorphism of the closed 2-dimensional annulus~$A$
twisting the two boundaries in opposite directions has at least two fixed points
in the interior of~$A$.
}

\begin{figure}[ht]
 \begin{center}
  \psfrag{A}{$A$}
  \psfrag{c}{$\chi (A) =0$}
 \leavevmode\epsfbox{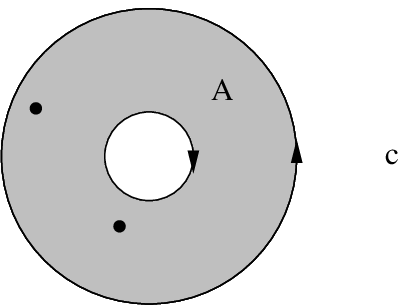}
 \end{center}
\end{figure}

\m  \ni
V.\ Arnold called this theorem ``the seed of symplectic topology''.
It has its origin in the planar circular restricted 3-body problem 
of celestial mechanics.
Poincar\'e, in his search for global periodic solutions, constructed a 
transversal section on the energy surface in the form of an annulus bounded
by the so-called direct and retrograde periodic orbits
and conjectured the above statement, that was later proved by
Birkhoff using two-dimensional methods.

The theorem is rather mysterious, and certainly not topological.
You cannot omit any of the assumptions.
To solve the mystery, Arnold suggested an alternative, more geometric
proof:
Assume that the maps in the theorem are, in addition, smooth and hence 
symplectic diffeomorphisms twisting the boundaries.
Take such a map~$\gf$, 
take two copies of the annulus endowed with this map, and
glue the annuli along their boundaries to obtain a 2-torus.
The map on the 2-torus defined by the two maps~$\gf$
can be smoothened near the glued boundaries to a symplectic diffeomorphism $\Phi$ that is Hamiltonian
and has no fixed points other than those of the two maps~$\gf$.
Now prove that $\Phi$ has at least three fixed points.
Then $\gf$ must have at least two fixed points.

\s
In order to explain the three fixed points in his conjecture, 
V.\ Arnold, familiar with the tricks of old mechanics, 
probably considered a symplectic diffeomorphism 
$\varphi \colon (x,y) \mapsto (X,Y)$ on~$T^2$ 
that is $C^2$-close to the identity, 
and hence, scrambling the variables, is represented by a single 
function~$G \colon T^2 \to \RR$ (the so-called generating function):
\begin{eqnarray*}
X &=& x + \frac{\partial G}{\partial Y}(x,Y), \\
y &=& Y + \frac{\partial G}{\partial x}(x,Y).
\end{eqnarray*}
The critical points of the function $G$ on $T^2$ are the 
fixed points of~$\varphi$ and hence by Ljusternik--Schnirelman,
this symplectic diffeomorphism on~$T^2$ possesses at least three fixed points.

\b \ni
{\bf 1982 -- Our proof \cite{CoZe83}}

\m \ni
C.\ Conley and I looked for forced oscillations of time-periodic
Hamiltonian systems. Our new idea was to use methods from dynamical systems
and to investigate the structure of the bounded orbits 
of the gradient flow of the action functional defined on the space of contractible
loops on the torus, because one expects (under suitable compactness conditions) 
that a bounded orbit of a gradient equation
converges automatically as time goes to plus infinity and minus infinity to critical points, 
which are the desired forced oscillations.
Since the gradient flow does not exist, we first reduced the problem to a variational problem
on a finite dimensional subspace of the loop space by the saddle point reduction method.

In our special case of the standard torus, it is convenient to first regularize the gradient and extend
the domain of definition $C^\infty(S^1, \RR^{2n})$ of the action functional~$f$ to the larger Sobolev space
$H^s (S^1, \RR^{2n})$ for $s = 1/2$.
We recall that the Sobolev spaces $H^s (S^1, \RR^{2n})$ for $s \geq 0$ are defined by
$$
H^s \,=\, \left\{ x \in L^2 (S^1,\RR^{2n}) \:\big|\, \sum_{k \in \ZZ}|k|^{2s} \, |x_k|^2 < \infty \right\}
$$
where 
\begin{equation} \label{e:Fourier}
x(t) \,=\, \sum_{k \in \ZZ} e^{k 2 \pi t J}x_k, \qquad x_k \in \RR^{2n}
\end{equation}
is the $L^2$-Fourier series of~$x$.
Then $H^s$ is a Hilbert space with inner product 
$$
\langle x,y \rangle_s \,:=\, \langle x_0,y_0 \rangle + 2\pi \sum_{k \in \ZZ}|k|^{2s} \langle x_k,y_k \rangle
$$
In the following we shall abbreviate $E = H^{1/2}$ and $\langle \cdot,\cdot \rangle = \langle \cdot,\cdot \rangle_{1/2}$
and $\| \; \| = \| \; \|_{1/2}$.
There is an orthogonal splitting $E = E^- \oplus E^0 \oplus E^+$ into the subspaces of~$E$
having only Fourier coefficients for $k<0$, $k=0$, and $k>0$, respectively.
So, every $x \in E$ has a unique decomposition
$$
x = x^- \oplus x^0 \oplus x^+ \,\in\, E^- \oplus E^0 \oplus E^+ .
$$
We identify the constant loops in $E^0$ with points of the torus $T^{2n}$.
Proceeding now as in \cite[Chapter~3]{HoZe94}, we first look at the symplectic part
of the functional~$f$, namely at
$$
a(x,y) \,:=\, \int_0^1 \frac 12 \langle -J \dot x, y \rangle \, dt, \qquad x,y \in C^\infty(S^1,\RR^{2n}) .
$$
Inserting the smooth loops $x(t)$ and $y(t)$ represented by their Fourier expansions~\eqref{e:Fourier}
and observing that
$$
\int_0^1 \langle e^{j 2\pi t J} x_j, e^{k 2\pi t J} x_k \rangle \,=\, \delta_{jk}\2 \langle x_j,x_k \rangle, 
$$
one computes that
$$
a(x,y) \,=\, \tfrac 12 \langle x^+, y^+ \rangle - \tfrac 12 \langle x^-, y^- \rangle
$$
for $x,y \in C^\infty (S^1, \RR^{2n})$. This expression is defined also for $x,y \in E$
and defines the extension of~$a$ to the Hilbert space~$E$.

Define the smooth function $a \colon E \to \RR$ by 
$$
a(x) = a(x,x) \,=\, \tfrac 12 \|x^+\|^2 - \tfrac 12 \|x^-\|^2 .
$$
The derivative of $a$ in $E$ is given by 
$da(x) \2 y = \langle x^+-x^-, y \rangle$, so that the $E$-gradient of~$a$ becomes
$$
\nabla a(x) \,=\, x^+-x^- \,\in\, E .
$$
For the regularized gradient of the Hamiltonian part 
$$
b(x) := \int_0^1 H(t,x(t)) \,dt, \quad x \in E
$$ 
one shows that
$$
\nabla b(x) \,=\, j^* \2 \nabla H(x) \in E, \quad \mbox{where } \nabla H (x) (t) := \nabla H (t, x(t)) \in L^2.
$$
Here, $j^* \colon L^2 \to E$ is the adjoint of the compact embedding $j \colon E \to L^2$, 
as usually defined by $\langle j(x), y \rangle_{L^2} = \langle x, j^*(y) \rangle_{1/2}$.
The functions $a,b \colon E \to \RR$ are smooth, and the regularized action functional
$$
f \colon E \to \RR, \quad f(x) = a(x)-b(x)
$$
is a smooth function.
Not all elements $x$ of $E = H^{1/2}$ are represented by continuous functions. 
However, if $x$ is a critical point of~$f$, i.e., if $x$ solves the equation
$$
\nabla f (x) \,=\, \nabla a(x) - \nabla b(x) \,=\, 0, 
$$
then $x \in C^\infty (S^1, \RR^{2n})$, and the loop $x(t)$ solves the Hamiltonian
equation and hence is a forced oscillation.

In order to find the critical points 
of $f \colon E \to \RR$ we take, following H.~Amann's saddle point reduction of~$f$,
the finite dimensional subspace $Z \subset E$ of smooth loops
$$
Z = E_{n_0} = \biggl\{ x \in E \mid x(t) = \sum_{|k| \leq n_0} e^{k 2\pi t J}x_k \biggr\}
$$
for a large integer $n_0$ to be determined.
Let $P = P_{n_0} \colon E \to Z$ be the orthogonal projection operator. 
The space~$E$ splits into $E = Z \oplus Y$, where $Z = P E$ and $Y = (1 - P)E$.
In order to solve $\nabla f(x) = 0$ for $x \in E$ we shall solve equivalently the 
two equations
\begin{equation} \label{e:P1P}
P \2 \nabla f(x) =0 \quad  \mbox{ and } \quad (1-P) \2 \nabla f(x) = 0
\end{equation}
for $x \in E$.
In order to solve the second equation of~\eqref{e:P1P} one first finds
for every $z \in Z$ a unique point $\varphi (z) \in Y$ solving the equation
\begin{equation} \label{e:z}
(1-P) \2 \nabla f(z+\varphi (z)) =0 
\end{equation}
provided that $n_0$ is sufficiently large.
For this one uses the compactness of~$j$ and~$j^*$ and the Banach fixed point theorem.
The function $\varphi \colon Z \to Y$ is smooth by the implicit function theorem,
and satisfies, by uniqueness, 
$\varphi (z+j) = \varphi (z)$ for all $j \in \ZZ^{2n}$ and all $z \in Z$.

\begin{figure}[ht]
 \begin{center}
  \psfrag{z}{$z$}
  \psfrag{Z}{$Z$}
  \psfrag{Y}{$Y$}
  \psfrag{f}{$\varphi (z)$}
 \leavevmode\epsfbox{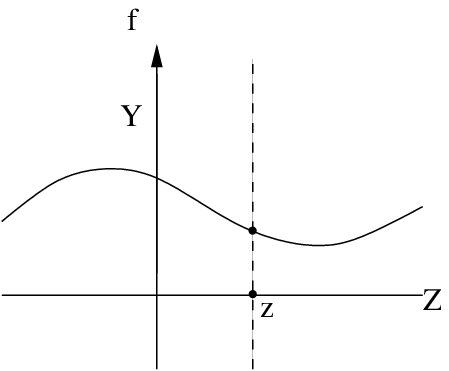}
 \end{center}
\end{figure}

\ni
Defining the smooth function $g \colon Z \to \RR$ by 
$$
g(z) := f(z+\varphi (z)) 
$$
one obtains, using~\eqref{e:z} and identifying the Hilbert spaces with their duals,
\begin{eqnarray*}
\nabla g(z) &=& (1 + d\varphi (z))^* \, \nabla f (z+\varphi (z)) \\
            &=& (1 + d\varphi (z))^* \, \left[ P \2 \nabla f + (1-P) \2 \nabla f \right] (z+\varphi (z)) \\
		&=& P \2 \nabla f (z+\varphi (z)) .
\end{eqnarray*}
To solve the first equation in~\eqref{e:P1P} it remains to find the critical points of the smooth function
$g \colon Z \to \RR$, that satisfies $g(z+j) = g(z)$ for all $j \in \ZZ^{2n}$ and $z \in Z$.
For this purpose, we look for all the bounded orbits of the gradient flow
$$
\frac{d}{ds} \psi^s (z) \,=\, \nabla g (\psi^s(z)) \quad \mbox{on $Z$.}
$$
Split $z = x \oplus \xi^+ \oplus \xi^- \in \left( E^0 \oplus E^+ \oplus E^- \right) \cap Z 
                          =: Z^0 \oplus Z^+ \oplus Z^-$,
and let $P^0, P^+, P^-$ be the orthogonal projections
of~$Z$ onto the spaces $Z^0, Z^+, Z^-$, respectively.
With respect to this splitting, write $\psi^s(z) = \left( x(s), \xi^+(s), \xi^-(s) \right)$.
Using the above identities for $\nabla g(z)$ and $\nabla a (z)$,
the flow $\psi^s$ is then represented by the ordinary differential equations on~$Z$
\begin{eqnarray*} 
\frac{d}{ds} x(s)\;\: &=& \phantom{-\xi^-(s)} -P^0 \2 \nabla b (\ast) ,\\
\frac{d}{ds} \xi^+(s) &=& \phantom{-}\xi^+(s)  - P^+ \2 \nabla b (\ast) ,\\
\frac{d}{ds} \xi^-(s) &=& -\xi^-(s) - P^- \2 \nabla b (\ast) ,
\end{eqnarray*}
where we abbreviated $\ast = \psi^s(z)+ \gf (\psi^s(z))$.
In view of the boundedness of~$\nabla b$, this immediately leads to the estimates
\begin{eqnarray*}
\frac{d}{ds} \2 \|\xi^+\|^2 &=& 2\, \Bigl\langle \frac{d}{ds} \xi^+, \xi^+ \Bigr\rangle  \,\geq 1 \quad \mbox{ if $\|\xi^+\| \geq K$},   \\
\frac{d}{ds} \2 \|\xi^-\|^2 &=& 2\, \Bigl\langle \frac{d}{ds} \xi^-, \xi^- \Bigr\rangle  \,\leq 1 \quad \mbox{ if $\|\xi^-\| \geq K$},   
\end{eqnarray*}
for a large constant $K>0$. We see that outside of the set
$$
B \,:=\, T^{2n} \times D^+ \times D^- \,\subset\, T^{2n} \times Z^+ \times Z^-,
$$
where $D^\pm \subset Z^\pm$ are the compact disks of radius~$K$, there are no bounded orbits!
By
$$
B^- = T^{2n} \times \pp D^+ \times D^-
$$
we denote the exit set of the flow in~$B$.
Schematically, the flow on~$Z$ looks as follows.

\begin{figure}[ht]
 \begin{center}
  \psfrag{S}{$S$}
  \psfrag{T}{$T^{2n}$}
  \psfrag{B}{$B$}
  \psfrag{x-}{$Z^-$}
  \psfrag{x+}{$Z^+$}
  \psfrag{dD}{$B^-$}
 \leavevmode\epsfbox{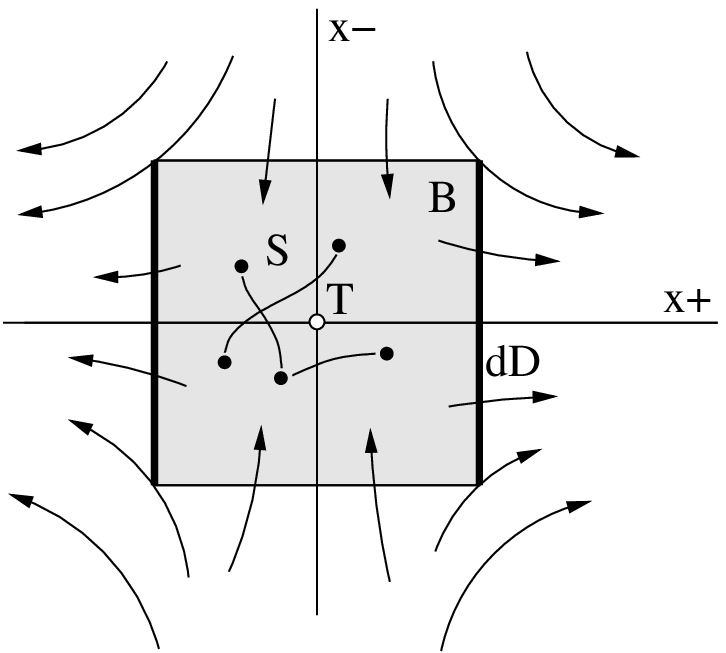}
 \end{center}
\end{figure}

The bounded orbits of the flow are all contained in the interior of $B$.
They constitute the maximal invariant set~$S$ of the compact space~$B$, 
formally defined by
$$
S \,=\, \left\{ z \in B \mid \varphi^s (z) \in B \mbox{ for all } s \in \RR \right\} .
$$
The set $S$ consists of the critical points of~$g$ together with the connecting orbits of the
flow of~$\nabla g$ between the critical points.
It is an invariant, compact, and isolated set in~$B$.

The compact pair $(B,B^-)$ is an example of a Conley index pair for the isolated
set~$S$. Collapsing the exit set~$B^-$ to a point~$\ast$, the homotopy type
$h(S)$ of~$S$, defined by 
$$
h(S) \,:=\, \left[ (B / B^-, \ast) \right] ,
$$
is called the Conley index of~$S$.
An algebraic invariant of~$h(S)$ is the cohomology of the compact pair $(B,B^-)$,
that is encoded in the Poincar\'e polynomial
$$
p(t,h(S)) \,=\, \sum_{j \geq 0} \dim \check H^j (B,B^-)\,t^j 
$$
where we use $\check{\rm C}$ech or Alexander--Spanier cohomology
with~$\ZZ_2$-coefficients.
We refer to \cite{Do80, ST94} for basic notions of algebraic topology.

The homotopy index is independent of the choice of the index pair
for the same invariant isolated compact set~$S$.
It is extremely stable under perturbations of the flow. 
In order to compute the index, Conley proved a crucial continuation
theorem for families of flows along which all the Conley indices
are homotopy equivalent. 
However, we are not interested in the homotopy type of~$h(S)$,
we are interested in the topology of the set~$S$ itself.
In order to describe it, I first recall two definitions.

\m \ni
{\bf Definition (Conley)}
A {\it Morse decomposition} of a compact invariant set~$S$
of a continuous flow on a metric space is an ordered finite 
collection $(M_j)_{j \in J}$ of disjoint compact and isolated invariant
subsets $M_j \subset S$ with the following property:
For every point
$$
p \in S \setminus \bigcup_J M_j
$$
there exists a pair $i<j$ of indices such that the positive limit set~$\omega (p)$
and the negative limit set $\omega^* (p)$ satisfy 
$$
\omega (p) \subset M_i \quad \mbox{and} \quad \omega^* (p) \subset M_j.
$$

Let me also recall the definition of the {\it cup-length} $\ell (X)$ of a non-empty
compact space~$X$:
$$
\ell (X) \,=\, 1 + \sup \left\{ k \in \NN \mid \exists \2 \ga_1, \dots, \ga_k  \mbox{ such that } \ga_1 \cup \dots \cup \ga_k \neq 0 \right\}
$$
where the $\ga_j \in \check H^*(X)$ must be graded elements of degree $\geq 1$.
If there is no such non-vanishing cup product, then $\ell (X) := 1$.

Using now some elementary results from algebraic topology, one derives the estimates
$$
2n+1 \,=\, \ell (T^{2n}) \,=\, \ell (B) \,\leq\, \ell (S) \,\leq\, \sum_{j \in J} \ell (M_j)
$$
for every continuous flow on $S$ having the index pair $(B,B^-)$ and for every Morse
decomposition of~$S$. 
In our case the flow is a gradient flow, and for the Arnold conjecture we can assume that it has 
only finitely many critical points~$z_j$.
Then the sets $M_j := \{z_j\}$ can be ordered such that they constitute a Morse decomposition of~$S$.
Since $\ell (\{z_j\}) =1$, we therefore conclude that
$$
\# \{\mbox{ critical points }\} \,\geq\, \ell (T^{2n}) \,=\, 2n+1 .
$$
We have found at least $2n+1$ critical points of $\nabla f$ and hence at least $2n+1$ forced oscillations.
If $n=1$, then $\ell (T^2) =3$ and the first part of the Arnold conjecture for the 2-torus is verified.

\s
In order to prove the second part of V.\ Arnold's conjecture for~$T^{2n}$, 
we assume that all the forced oscillations of the Hamiltonian system are non-degenerate. 
Then all the critical points of~$g$ are non-degenerate and finite in number.
We can then apply the Morse inequalities from \cite[Theorem 3.3]{CoZe84}
for the cohomology of an ordered Morse decomposition $(M_j)$ of~$S$. They are given by
$$
\sum_{j \in J} p(t,h(M_j)) \,=\, p (t,h(S)) + (1+t) Q(t) 
$$
where the polynomial $Q$ has only non-negative integer coefficients.
The proof is straightforward and follows from elementary dynamical systems methods;
it does not require a manifold, but just a continuous flow
(not necessarily a gradient flow) on a compact metric space~$S$ and a Morse decomposition of this flow.

Recalling that $(B,B^-)$ is an index pair of~$S$, we can apply the K\"unneth formula 
and obtain
\begin{eqnarray*}
\check H^* (h(S)) \,=\, 
\check H^* (B,B^-) &=& 
\check H^* (T^{2n} \times D^+, T^{2n} \times \pp D^+) \\
&=& \check H^* (T^{2n}) \otimes  \check H^* ( D^+, \pp D^+) \\
&=& \check H^* (T^{2n}) \otimes  \check H^* (\dot S^N) .
\end{eqnarray*}
Here, $\dot S^N$ denotes a sphere of dimension $N = \dim Z^+$
with a distinguished point $\ast$, that is, $\dot S^N = (S^N,\ast)$. 
Therefore, 
$$
p(t,h(S)) \,=\, \sum_{j =0}^{2n} \binom{2n}{j} \, t^j \, t^N.
$$
The finitely many non-degenerate critical points $\{z_j\} = M_j$ of~$g$ on~$S$ 
can be ordered to constitute a Morse decomposition of~$S$.
The critical points are hyperbolic and hence 
$$
p(t,h(z_j)) \,=\, p(t,\dot S^{m_j}) \,=\, t^{m_j}
$$  
where $m_j$ is the Morse index of~$z_j$. 
To sum up, the Morse inequalities in our special case are
$$
\sum_{j} t^{m_j} \,=\, \sum_{j=0}^{2n}  \binom{2n}{j} \, t^{j+N} + (1+t)Q(t).
$$
We read off that the number of non-degenerate critical points of~$g$
is at least 
$$
\sum_{j=0}^{2n}  \binom{2n}{j} \,=\,2^{2n}
$$
which is the sum of the Betti numbers of $T^{2n}$.
If $n=1$, then $2^2=4$, the lower bound conjectured by V.\ Arnold for the 2-torus.
This ends the sketch of our proof.

\s
One sees that the structure of the bounded orbits $S$ of the gradient flow automatically
leads to the Arnold conjecture on the torus.
The forced oscillations found are related to the topology of the underlying manifold, 
not to the topology of the loop space as in the geometric problem of 
closed geodesics on a compact Riemannian manifold.
This is due to the symplectic structure.

\s
After a talk about this rather simple result, we received
an unexpected compliment from Misha Gromov. 
He approached us and said ``This result is so beautiful, 
I would not mind stealing it!''

\s
While there are by now simpler proofs of the Arnold conjecture in the special case of the torus,
see for instance~\cite{Chap84},
the above proof ultimately lead to Floer homology and to the proof of the Arnold conjecture
for general symplectic manifolds by Andreas Floer.

\b \ni
{\bf Immediate consequences}

\s \ni
Right after the appearance of the proof of the fixed point theorem for $T^{2n}$,
Marc Chaperon in his Bourbaki Seminar talk 1982--83 observed that the proof also allows 
to verify another longstanding conjecture of V.\ Arnold about the number of 
intersections of Lagrangian submanifolds in the special case of a Hamiltonian map
$\varphi^1$ on the symplectic manifold $T^* T^m$, namely
$$
\# \bigl(T^m \cap \varphi^1(T^m) \bigr) \,\geq\, m+1
$$
resp.\ $\geq 2^m$ if all the intersections are transverse, see \cite{Chap83}.

Moreover, at the same time Michel Herman used the fixed point result for~$T^{2n}$
to show in a letter to me, see also~\cite{Herman}, that the set of smooth symplectic 
diffeomorphisms on~$T^{2n}$ in the set of volume preserving diffeomorphisms,
$$
C^\infty_{\symp}(T^{2n},\omega_0) \,\subset\, C^\infty_{\vol} (T^{2n},\omega_0) ,
$$ 
is not $C^0$-dense for $n>1$.
Thus he discovered that the set of symplectic diffeomorphisms and the set of volume preserving 
diffeomorphisms are quite different, if $n>1$.

\s
As you all know, all these pioneering results have been extended later on to more general 
symplectic manifolds.
Of course, with Charley Conley we did try do apply our ideas to general compact symplectic
manifolds. 
We wrote down on $\RR^{2n}$ the equation of a bounded orbit of the not regularized $L^2$ gradient flow
of the action functional
and ended up with a non-linear asymptotic boundary value problem for an elliptic system of PDEs, 
but did not know how to proceed from there.

\s
\begin{figure}[ht]
 \begin{center}
  \includegraphics[width=4.5cm]{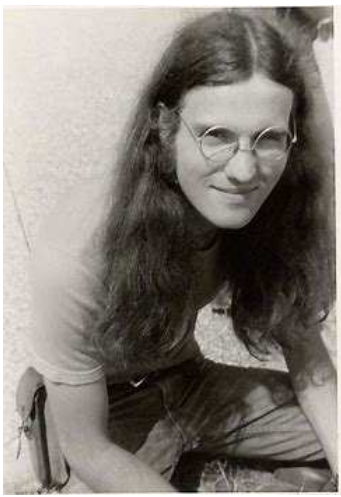}
 \end{center}
 \caption{Andreas Floer in 1976}
\end{figure}


\s
One day a student knocked on my office door on the 7th floor.
I had never seen him in my Analysis lectures.
He did not introduce himself but just asked: 
``Do you have an interesting topic for my thesis?''
I did not know what he meant by ``interesting'' and so 
asked him in, and explained in detail the Arnold conjecture, 
its background, our proof for the special case of the torus, 
the idea to analyze the structure of bounded solutions of the
artificial gradient flow, and also the use of Conley's index theory.
He was immediately attracted to the challenging problem, 
and I suggested that he should start with the Arnold conjecture
for surfaces~$\Sigma_g$ of higher genus~$g \geq 2$ in order to verify the claims of Eliashberg.
And that is exactly what he did. Very soon, he showed up in my office
with beautiful new ideas. Only his notes were quite a mess.
So I told him to write his ideas down with loving care, 
for the benefit of the reader.
At first, he was furious, but later on he always made fun of me with 
``my loving care writing style''.
But then he did it, after finitely many iterations.
Whenever he showed up with a new correction, he said, before even entering my office:
``Herr Zehnder, Sie k\"onnen mich jetzt wieder in der Luft zerreissen''
(now you can tear me apart again).

In his proof of the Arnold conjecture for surfaces in~\cite{Flo86:surface},
Floer also reduced the variational functional, from an $L^2$-bundle
to a finite dimensional sub-bundle, and substantially refined the Conley index theory
by adding a new structure to the Conley index, 
that he used to prove a continuation theorem for the invariant set itself (not only for the index pair), 
see~\cite{Floer87}.
This allowed him to verify that the cohomology of~$\Sigma_g$ injects into the cohomology of the set of bounded orbits.
From this one finds the conjectured
estimates for the number of forced oscillations on~$\Sigma_G$ like for the torus, 
namely at least as many as the cup-length of the surface, which is three, 
and in the non-degenerate case at least as many as the sum of the Betti numbers, 
which is $2+2g$.

This convinced me that Andreas has the vision and the power to solve the general Arnold conjecture,
and I urged him to do so, what, among many other striking things, 
he did in a spectacular and revolutionary manner, 
extending Conley's ideas in his own way.
He had to explain me his approach several times and always started his explanations with:
``Now, Herr Zehnder, I explain it to you for the $n$th time.''

\s
The only person who did not believe in Floer's ideas was M.\ Gromov,
because they did not immediately fit into his framework of mathematical thinking. 
Typically for Andreas, such doubts did not make him angry or insecure, 
he simply laughed about them.

\m
Let me quickly recall some of Andreas' ideas, omitting
all the details that are technically very intricate,
involving the glueing techniques that he learned from Cliff Taubes
in Berkeley as well as elliptic PDE techniques, non-trivial non-linear functional analysis,
and subtle compactness considerations.

\b \ni
{\bf Floer's proof of the V.\ Arnold conjecture in the non-degenerate case~\cite{Flo89:hol}}

\m \ni
For a smooth time-periodic Hamiltonian function $H \colon S^1 \times M \to \RR$
on a compact symplectic manifold $(M,\omega)$,
Arnold conjectured in 1972, see \cite{Arn72} and also~\cite{Arn78}, 
that 
$$
\# \left\{\mbox{ fixed points of $\varphi_H^1$ on $M$ }\right\} \,\geq\, 
\mbox{sum of the Betti numbers of $M$}
$$
under the assumption that all the contractible forced oscillations are non-degenerate.

In contrast to the torus case, in his proof of this conjecture 
for so-called monotone symplectic manifolds, 
Andreas Floer avoided any finite dimensional reduction and started directly 
with the structure of all bounded orbits of the (not regularized) gradient of the action functional~$f$, 
defined on the set $\Omega$ of smooth contractible parametrized loops
$u \colon S^1 \to M$.
A gradient flow for this functional does not exist.
But the structure of all bounded orbits lead him to his Floer homology, 
whose chain complex is generated by the finitely many forced oscillations.
Recalling the continuation theorem for the Conley index, Floer proved in a second step that
his homology is independent of the choice of the Hamiltonian function~$H$.
Finally, he showed that for special Hamiltonians, namely time-independent and $C^2$-small Morse functions,
the Floer homology is the Morse--Smale homology of the manifold~$M$ and hence isomorphic to $H_*(M)$.

In our more detailed sketch of Floer's proof we require for simplicity 
that the first Chern class $c_1$ and the cohomology class $[\omega]$ 
vanish on the second homotopy group $\pi_2(M)$, 
and follow the exposition in~\cite{SaZe92}
and~\cite[Chapter~6, Section~5]{HoZe94}.

The time-dependent Hamiltonian function $H(t,x) = H(t+1,x)$ defines the time-dependent Hamiltonian
vector field $X_H$ by 
$$
\omega (X_H, \cdot) = -dH(\cdot) .
$$
We choose an almost complex structure~$J$ on~$M$ (namely an endomorphism of the tangent bundle such that $J^2 = -1$)
with the property that
$$
\omega_x (\xi, J(x)\eta) = g_x(\xi, \eta), \qquad \xi,\eta \in T_xM
$$
defines a Riemannian metric $g = g_J$ on~$M$.
Then the Hamiltonian vector field is represented as
$$
X_H(t,x) \,=\, J(x) \2 \nabla H(t,x), \qquad x \in M ,
$$
where the gradient is taken with respect to the metric~$g_J$.
The $L^2$-gradient of the action functional~$f$ on the space~$\Omega$ of contractible loops~$u$ in~$M$ 
then becomes
$$
\nabla f (u) \,=\, J(u) \2 \dot u + \nabla H(t,u), \qquad u(t) = u(t+1) .
$$
A bounded orbit of the gradient flow is a solution of the equation on~$\Omega$
$$
\frac{\pp}{\pp s} u(s) \,=\, -\nabla f(u(s)), \qquad s \in \RR
$$
where $u(s)(t) := u(s,t) = u(s, t+1)$
has bounded energy
$$
E(u) \,:=\, \frac 12 
\int_\RR 
 \int_{S^1} \left( \left\| \frac{\pp u}{\pp s} \right\|^2 + \left\| \frac{\pp u}{\pp t} - X_H(t,u) \right\|^2 \right) \, dt \2 ds  
\,<\, \infty .
$$
Here, the norms are induced by the Riemannian metric $g_J$.
Explicitely, the function $u \colon \RR \times S^1 \to M$ is a solution of the elliptic partial differential equation
$$
\frac{\pp u}{\pp s} + J(u) \2 \frac{\pp u}{\pp t} + \nabla H (t,u) \,=\, 0 .
$$
The requirement of bounded energy of the solution~$u$ is, as expected, equivalent to the asymptotic boundary
conditions
$$
\lim_{s \to -\infty} u(s,t) = x_-(t) \quad \mbox{ and } \quad \lim_{s \to +\infty} u(s,t) = x_+(t)
$$
for two non-degenerate forced oscillations $x_-$ and $x_+$.
We obtain a so-called Floer connecting orbit of contractible loops:

\begin{figure}[ht]
 \begin{center}
  \psfrag{s}{$s$}
  \psfrag{x-}{$x_-$}
  \psfrag{x+}{$x_+$}
 \leavevmode\epsfbox{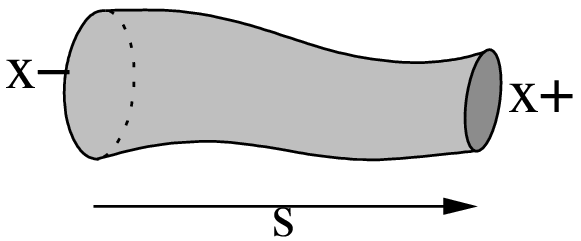}
 \end{center}
\end{figure}
%

\ni
We are now confronted with a non-linear elliptic system of PDEs for $u(s,t)$
with asymptotic (in $s$) and periodic (in $t$) boundary conditions.
Linearizing the differential operator along a solution~$u$ we arrive at 
the Fredholm operator
$$
F(u) \,=\, \nabla_{\!s}\2 \xi + J(u) \2 \nabla_{\!t}\2 \xi + \nabla_{\!\xi} \2 (J (u)) \2 \frac{\pp u}{\pp t} + \nabla_{\!\xi} \nabla H(t,u)
$$
between appropriate Sobolev spaces, where $\xi \in C^\infty (u^* TM)$
and the covariant derivatives are defined by the metric~$g_J$.
For a good choice of~$J$ the linear Fredholm operator~$F(u)$ is surjective by the Sard--Smale theorem
(in the proper setting).
Hence there exists a smooth finite dimensional manifold $\cm (x_-,x_+)$ of Floer connecting orbits
whose dimension is equal to
$$
\dim \cm (x_-,x_+) \,=\, \Fred (u) \,=\, \mu (x_-) - \mu (x_+) .
$$
Further, by compactness, all these manifolds have  finitely many components.

Next, we introduce the finite dimensional graded $\ZZ_2$-vector space~$C$
generated by the non-degenerate forced oscillations $\cp_H$:
$$
C \,=\, \oplus C_k
$$
where $C_k = \span_{\ZZ_2} \{ x \in \cp_H \mid \mu (x) = k \}$.
On this vector space Floer defined a boundary operator $\pp \colon C \to C$
by defining it on generators $x \in C_k$ by
$$
\pp_k (x) \,=\, \sum_{\mu (y) = k-1} \langle \pp x, y \rangle \, y 
$$
where $\langle \pp x, y \rangle$ is the parity of the finite number of components
of the 1-dimensional manifold $\cm (x,y)$.

Using the glueing technique that he had learned from Cliff Taubes 
and an implicit function theorem in Banach spaces, 
Floer verified that $\pp \circ \pp =0$, i.e., $\pp_{k-1} \circ \pp_k =0$ for all~$k$,
or explicitely, 
$$
\sum_{\mu (z) = k-2} \left( \sum_{\mu (y) = k-1} \langle \pp x, y \rangle \langle \pp y, z \rangle \right) z \,=\, 0 .
$$
To see that for each pair $(x,z)$ the sum in the large bracket vanishes mod~2, 
Floer observed that the 1-dimensional manifold $\cm (x,z) / \RR$, 
obtained by taking the quotient of the 2-dimensional manifold~$\cm (x,z)$ by the free shift-action in~$s$,
has an even number of ends,  
as illustrated in the following picture.

\begin{figure}[ht]
 \begin{center}
  \psfrag{x}{$x$}
  \psfrag{y}{$y$}
  \psfrag{y'}{$y'$}
  \psfrag{z}{$z$}
  \psfrag{u}{$u$}
  \psfrag{u'}{$u'$}
  \psfrag{v}{$v$}
  \psfrag{v'}{$v'$}
  \psfrag{mu}{$\cm(x,z)$:}
  \psfrag{C1}{$C_{k}$}
  \psfrag{C2}{$C_{k-1}$}
  \psfrag{C3}{$C_{k-2}$}
 \leavevmode\epsfbox{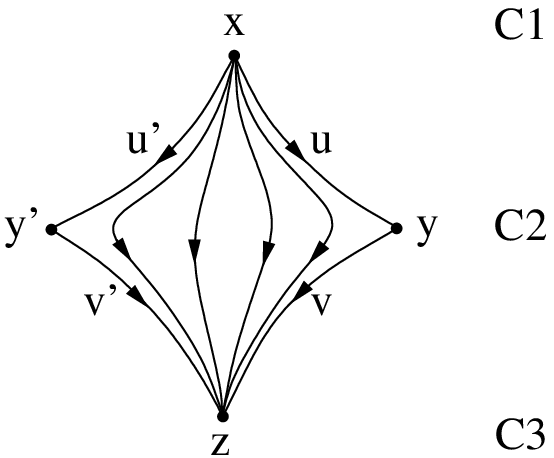}
 \end{center}
\end{figure}
%

\ni
The relevant broken connecting orbits at the ends therefore appear in pairs and so
$\partial \circ \partial = 0$ mod~2. 
The homology of the complex $(C,\partial)$,
$$
\FH_k (M,H,J) \,=\, \frac{\ker (\pp \colon C_k \to C_{k-1})}{\im (\pp \colon C_{k+1}\to C_k)} ,
$$
is called Floer homology.
Floer himself called it in~\cite{Flo89:hol} the {\it homological Conley index}\/
of the invariant set of bounded orbits of the gradient flow of the action functional.

Inspired by Conley's continuation theorem for the Conley index, 
Floer showed in a second step that his Floer homology is independent of the choice 
of~$H$ and~$J$. He takes two Hamiltonians $(H^\ga, J^\ga, x^\ga)$ and $(H^\gb, J^\gb, x^\gb)$
with the associated almost complex structures and generators of the corresponding 
chain complexes, and defines a clever homotopy between the Hamiltonians and almost complex structures,
which satisfies, in particular, 
\begin{equation*}
H(s,t,x) \,=\, \left\{ 
\begin{array}{cc}
H^\ga (t,x) & \mbox{as } s \to -\infty  \\ [0.2em]
H^\gb (t,x) & \mbox{as } s \to +\infty 
\end{array}
\right.
\end{equation*}
and similarly for~$J$.
Then he studies the $s$-dependent solutions of the PDE
$$
\frac{\pp u}{\pp s} + J(s,u) \2 \frac{\pp u}{\pp t} + \nabla H(s,t,u) = 0
$$
with asymptotic boundary conditions
$$
\lim_{s \to -\infty} u(s,t) = x^\ga(t) \quad \mbox{ and } \quad
\lim_{s \to +\infty} u(s,t) = x^\gb(t) .
$$
The manifold $\cm (x^\ga, x^\gb)$ of Floer connecting orbits has, this time, 
dimension
$$
\dim \cm (x^\ga, x^\gb) \,=\, 
         \mu (x^\ga, H^\ga) - \mu (x^\gb, H^\gb) - 1.
$$
Next, Floer constructs a homomorphism of the complexes
$$
\phi^{\gb \ga}_k \colon C_k(M,H^\ga) \to C_k(M,H^\gb)
$$
satisfying $\pp^\gb \circ \phi^{\gb \ga} = \phi^{\gb \ga} \circ \pp^\ga$.
It respects the grading and induces an isomorphism between the Floer homologies
$$
\Phi^{\gb \ga}_k \colon \FH_k(M,H^\ga,J^\ga) \to \FH_k(M,H^\gb,J^\gb)
$$
satisfying $\Phi_*^{\gg \gb} \circ \Phi_*^{\gb \ga} = \Phi_*^{\gg \ga}$ and $\Phi_*^{\ga \ga} = \id$.

Finally, in order to show that the Floer homology is useful, 
Floer computes it. For this purpose, he takes a time-independent Hamiltonian function 
$h \colon M \to \RR$ which is a $C^2$-small Morse function.
Its forced oscillations are independent of~$t$ and are the critical points of~$h$ on~$M$.
Moreover, the Floer connecting orbits $u(s,t)$ are also independent of~$t$,
hence $u(s,t) = \gamma (s)$.
It follows that the PDE for bounded orbits of~$f$ in~$\Omega$ 
reduces to the ODE of the gradient equation of the Morse function~$h$
$$
\frac{d}{ds} \gamma (s) \,=\, - \nabla h (\gamma (s)) 
$$
on the underlying compact symplectic manifold~$M$.
The $\mu$-index of a critical point~$x$ of~$h$ as a forced oscillation 
is related to the Morse index of~$x$ by
$$
\mu (x,h) \,=\, \ind_h(x) - n 
$$
where $2n = \dim M$.
Since, in addition, the Floer boundary operator is equal to the Morse--Smale boundary operator, we conclude that
\begin{eqnarray*}
\FH_*(M,h,J) &\cong& \mbox{Morse--Smale homology of $h$ on $M$} \\
 &\cong& \mbox{singular homology of $M$.}
\end{eqnarray*}
Therefore, and in view of the independence of the Hamiltonian, 
we conclude that for every $H$ on~$M$,
$$
\#  \left\{\mbox{ fixed points of $\varphi_H^1$ on $M$ }\right\} \,\geq\, 
\mbox{sum of the Betti numbers of $M$}.
$$
This finishes the sketch of Andreas Floer's proof of the general 
Arnold conjecture in the non-degenerate case.

\s
\begin{figure}[ht]
 \begin{center}
  \includegraphics[width=5cm]{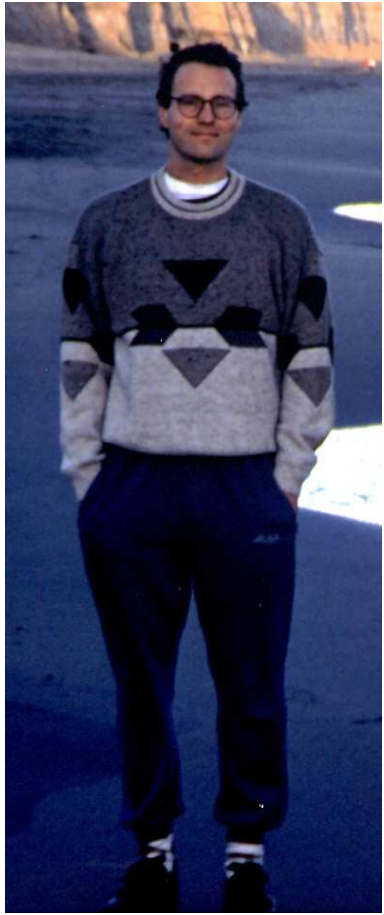}
 \end{center}
 \caption{Andreas Floer in 1989}
\end{figure}
%

\m
Andreas Floer's methods, ideas, and constructions were a crucial break-through and continue
to influence symplectic topology and Hamiltonian dynamics enormously, 
see for instance the subsequent survey~\cite{AbSch19}.

Over the years I have followed with great interest and joy the dynamic development
of the field of symplectic topology at the RUB and the creation of the Floer centre.
I stop here with my memories of the beginnings of symplectic topology more
than thirty years ago, here in Bochum.


\end{document}